\documentclass[twoside,10pt]{amsart}
\usepackage{amsmath,amsthm,amsfonts,graphicx,listings,centernot,fullpage}
\usepackage[english]{babel}
\usepackage{latexsym}
\usepackage{amssymb}
\usepackage{color}
\usepackage{graphicx}

\newtheorem{theorem}{Theorem}

\newtheorem{problem}{Problem}
\newtheorem{conjecture}{Conjecture}

\newtheorem{lemma}[theorem]{Lemma}
\newtheorem{corollary}[theorem]{Corollary}
\newtheorem{proposition}[theorem]{Proposition}

\def\Proof{\textbf{Proof.  }}
\newcommand{\open}{{\mathcal O}}
\newcommand{\dense}{{\mathcal D}}

\newcommand{\sone}{{\sf S}_1}
\newcommand{\sfin}{{\sf S}_{fin}}
\newcommand{\gone}{{\sf G}_1}

\newcommand{\naturals}{{\mathbb N}}
\newcommand{\reals}{{\mathbb R}}
\newcommand{\rationals}{{\mathbb Q}}

\newcommand{\sorgenfrey}{{\mathbb S}}

\newcommand{\cohen}{{\mathbb C}}

\address{Department of Mathematics\\ Boise State University\\ Boise, Idaho 83725}
\email{liljanababinkostova@boisestate.edu,\\
       mscheepe@boisestate.edu}
\address{Dipartimento di Matematica, Universita di Messina, Via F. Stagno d'Alcontres N.31,
          8166 Messina (Italy). (Pansera) }
\email{bpansera@unime.it}

\begin{document}

\title{Weak covering properties and selection principles}

\author{L. Babinkostova, B. A. Pansera and M. Scheepers}

\date{}

\maketitle

\begin{abstract} No convenient internal characterization of spaces that are productively Lindel\"of is known. Perhaps the best general result known is Alster's internal characterization, under the Continuum Hypothesis, of productively Lindel\"of spaces which have a basis of cardinality at most $\aleph_1$. It turns out that topological spaces having Alster's property are also productively weakly Lindel\"of. The weakly Lindel\"of spaces form a much larger class of spaces than the Lindel\"of spaces. In many instances spaces having Alster's property satisfy a seemingly stronger version of Alster's property and consequently are productively X, where X is a covering property stronger than the Lindel\"of property. This paper examines the question: When is it the case that a space that is productively X is also productively Y, where X and Y are covering properties related to the Lindel\"of property.
\end{abstract}

\section{Introduction}

A topological space is said to be Lindel\"of if each of its open covers contains a countable subset that covers the space. Though this class of spaces has been extensively studied there are still several easy to state problems that have not been resolved. Call a space $X$ \emph{productively Lindel\"of} if $X\times Y$ is a Lindel\"of space whenever $Y$ is a Lindel\"of space. 

In the quest to find an internal characterization of the productively Lindel\"of spaces K. Alster identified the following conditions: Call a family $\mathcal{F}$ of ${\sf G}_{\delta}$ subsets of a space $X$ a ${\sf G}_{\delta}$ \emph{compact} cover if there is for each compact subset $K$ of $X$ a set $F\in\mathcal{F}$ such that $K\subseteq F$. A space is said to be an \emph{Alster space} if each ${\sf G}_{\delta}$ compact cover of the space has a countable subset covering the space\footnote{This definition is not identical to the definition in \cite{BKR}, but is equivalent to it, and is in fact the property (*) defined by Alster on p. 133 of \cite{Alster1}.}. Alster (and independently \cite{BKR}) proved that if $X$ is an Alster space then it is productively Lindel\"of, and $X^{\aleph_0}$ is Lindel\"of. 
\begin{problem}[K. Alster]\label{Alsterpr} Is every productively Lindel\"of space an Alster space?
\end{problem}

A significant body of partial results has developed around this problem, yet no definitive answer is known to it. Several weakenings of the Lindel\"of property that have been investigated because of their natural occurrence in some mathematical contexts. The corresponding product theory for these is not as extensively developed as for Lindel\"of spaces. Product theoretic questions may be more manageable for the corresponding weakened analogues of the selective versions of the Lindel\"of property.

In another direction, a number of selective versions of the Lindel\"of property and weakenings of it have been investigated because of their relevance to several other mathematical problems. It has been found that soome questions about Lindel\"of spaces are ``easier" for these more restricted classes. It is natural to inquire whether the corresponding product-theoretic problem for these narrower classes of spaces is more manageable. Some such questions have been raised: For example: In \cite{TallNote} and \cite{TallTsaban} the notion of a \emph{productively Menger} space is considered and in \cite{AT} the notion of a \emph{productively FC-Lindel\"of space} is introduced.

And thirdly, solving versions of a problem still unresolved for Lindel\"of spaces by strengthening the hypotheses to selective versions while weakening the conclusions to weak covering properties, may yield some insights on the original problem.
Progress on the internal characterization problem may also yield insights on an older problem of E. Michael:
\begin{problem}[E.A. Michael]\label{prodLind} If $X$ is a productively Lindel\"of space, then is $X^{\aleph_0}$ a Lindel\"of space?
\end{problem}

These are the motivations for our paper. The paper is organized as follows. In Section 2 we introduce some basic notation and terminology. In Section 3 we give a number of examples of when products fail to have some of the properties we are investigating. Section 4 focuses on the question of characterizing the productively Lindel\"of spaces. We raise the question of when {\sf P} and {\sf Q} are covering properties of topological spaces, is a space that is productively {\sf P} also productively {\sf Q}? In Section 5 we incorporate the weak covering properties into the investigation.


\section{Some notation and terminology}

A space is said to be \emph{weakly Lindel\"of} if each of its open covers contains a countable subset for which the union is dense in the space. A space is said to be \emph{almost Lindel\"of} if each of its open covers contains a countable subset for which the set of closures of elements of the countable set is a cover of the space. Both of these properties are weaker than the Lindel\"of property, and we have the following implications:
Lindel\"of $\Rightarrow$ almost Lindel\"of $\Rightarrow$ weakly Lindel\"of. For spaces with the ${\sf T}_3$ separation property almost Lindel\"of also implies Lindel\"of. Aside from this there are no other implications among these three properties. We will focus on the ``weakly" properties in this paper, leaving the ``almost" properties for another time.

Next we describe selective versions of these covering properties. We use the following notation for three of several upcoming relevant classes of families of open sets of a given topological space:\\

\begin{center}
\begin{enumerate}
\item[$\open$]\hspace{0.5in}The collection of open covers
\item[$\dense$]\hspace{0.45in} $\{\mathcal{U}: (\forall U\in\mathcal{U})(U$ open) and $(\bigcup\mathcal{U}$ dense in $X)\}$
\end{enumerate}
\end{center}

Let $\mathcal A$ and $\mathcal B$ be collections of subsets of an infinite set. Then
$\sone({\mathcal A}, {\mathcal B})$ denotes the following hypothesis:
\begin{quote}
For each sequence $(A_n : n \in {\mathbb N})$ of elements of $\mathcal A$ there is a sequence $(B_n : n\in {\mathbb N})$ such that, for each $n$, $B_n\in A_n$ and $\{B_n : n \in {\mathbb N}\}$ is an element of $\mathcal B$.
\end{quote}
Thus, $\sone(\open,\open)$ denotes the classical \emph{Rothberger} property. We shall call spaces with the property 
$\sone(\open,\dense)$ \emph{weakly Rothberger}.

The symbol $\sfin({\mathcal A}, {\mathcal B})$ denotes the hypothesis
\begin{quote}
For each sequence $(A_n : n \in {\mathbb N})$ of elements of $\mathcal A$ there is a sequence $(B_n : n\in {\mathbb N})$ such that, for each $n$, $B_n\subseteq A_n$ is finite, and $\bigcup\{B_n : n \in {\mathbb N}\}$ is an element of $\mathcal B$.
\end{quote}
$\sfin(\open,\open)$ denotes the classical \emph{Menger} property, while 
$\sfin(\open,\dense)$ denotes the \emph{weakly Menger} property.


















Several additional families $\mathcal{A}$ and $\mathcal{B}$ of topologically significant objects will be introduced as needed during the rest of the paper.
Our conventions for the rest of the paper are: By ``space" we mean a topological space. Unless other separation axioms are indicated specifically, we assume all spaces to be infinite and ${\sf T}_1$. Undefined notation and terminology will be as in \cite{E}.

\section{Possibilities of failure for products}

Towards investigating the product theory as outlined above, we consider if it is possible for certain products to fail having a covering property. We asked above, for example, if there could be a Rothberger space whose product with the space of irrational numbers is not weakly Lindel\"of. First, we settle that it is at least possible that the product of two Rothberger spaces can fail to be a weakly Lindel\"of space. We give two examples of how this could be. Both are consistency results.

\begin{theorem}\label{rothbweakl} It is consistent, relative to the consistency of {\sf ZFC}, that there are Rothberger spaces $X$ and $Y$ for which $X\times Y$ is not weakly Lindel\"of.
\end{theorem}
\Proof In \cite{HJ} Hajnal and Juhasz give examples $X$ and $Y$ of Lindel\"of spaces for which $X\times Y$ is not weakly Lindel\"of\footnote{A nice exposition of this example can be found in \cite{Stoyanova}, Example 3.25.}. These examples are constructed in {\sf ZFC}. Now consider these two ground model examples in the generic extension obtained by adding $\kappa>\aleph_0$ Cohen reals. Since $X$ and $Y$ are Lindel\"of in the ground model, Theorem 11 of \cite{MSFT} implies that $X$ and $Y$ are Rothberger in the generic extension. Since $X\times Y$ is not weakly Lindel\"of in the ground model, Theorem 1 of \cite{BaPaSc} implies that $X\times Y$ is not weakly Lindel\"of in the generic extension.
$\Box$

Our second example is a little stronger than the one just given. A \emph{Souslin line} is a complete dense linearly ordered space $X$ which is not separable but every family of disjoint intervals is countable. {\sf SH}, the \emph{Souslin Hypothesis}, states that there are no Souslin lines. {\sf SH} is independent of {\sf ZFC}. Theorem \ref{rothbweakl} is proven using forcing. One may ask to  what degree axiomatic circumstances determine whether a product fails to have the covering property of its factor spaces. Shelah \cite{Shelah84} proved that in the generic extension obtained by adding a Cohen real there is a Souslin line. Thus, the following Theorem \ref{notSHth} improves Theorem \ref{rothbweakl}. In the proof of this theorem we use the following notation. If $(L,<)$ is a linearly ordered set there are three topologies considered on it: We denote the topological space by $L$ if the topology is generated by sets of the form $(a,b)$ where $a<b$ are elements of $L$. When the topology is generated instead by sets of the form $\lbrack a,\, b)$, the topological space is denoted by the symbol $L^+$. When the topology is generated by sets of the form $(a,b\rbrack$, then the topological space is denoted by the symbol $L^-$.
\begin{theorem}[$\neg${\sf SH}]\label{notSHth} There are Rothberger spaces $X$ and $Y$ such that $X\times Y$ is not weakly Lindel\"of.
\end{theorem}
\Proof Let $L$ be a Souslin line. We may assume that $L$ has no nonempty open intervals that are separable. Then $L^+$ as well as $L^-$ are Lindel\"of spaces (\cite{Stoyanova}, Lemma 3.31). But $L^+$ is a refinement of the standard topology on $L$, and so $L$ is Lindel\"of.

{\flushleft{\bf Claim 1: }} $L^+$ (and similarly $L^-$) is a Rothberger space\footnote{Towards proving this we refine the argument from page 19 of \cite{COC5}.}.

Let $(\mathcal{U}_n:n\in\naturals)$ be a sequence of open covers of $L^+$. By Lemma 3.31 in \cite{Stoyanova} we may assume that each $\mathcal{U}_n$ consists of countably many open intervals of form $\lbrack a^n_k,b^n_k)$, $k\in\naturals$.

Consider $(\mathcal{U}_{2\cdot n}:n\in\naturals)$. The set $A=\overline{\{a^{2\cdot n}_k:k\in\naturals\}\bigcup\{b^{2\cdot n}_k:k\in\naturals\}}$ is countable (lest $L$ has a separable uncountable interval) nowhere dense in $L$. Choose for each $x\in L\setminus A$ an open interval $I_x\subseteq L\setminus A$ with $x\in I_x$. Note that for each $n$ there is a $\kappa$ with $I_x\subseteq (a^{2\cdot n}_\kappa,\, b^{2\cdot n}_\kappa)$ (since $I_x\cap A = \emptyset$).

Choose from each $\mathcal{U}_{2\cdot n+1}$ an element $J_{2\cdot n+1}$ such that $A\subseteq \bigcup_{n\in\naturals}J_{2\cdot n+1}$. Now $\kappa= L\setminus (\bigcup_{n\in\naturals}J_{2\cdot n+1})$ is a closed subset of $L$ and so Lindel\"of. Moreover $\{I_x: x\in \kappa\}$ is an open cover of $\kappa$, and so has a countable subcover, say $\{I_{x_n}:n\in\naturals\}$. Now choose for each $n$ a $\kappa_n$ such that $I_{x_n}\subseteq J_{2\cdot n} = \lbrack a^{2\cdot n}_{\kappa_n},\, b^{2\cdot n}_{\kappa_n}) \in \mathcal{U}_{2\cdot n}$.
Then the sequence $(J_n:n\in\naturals)$ is a cover of $L^+$, and for each $n$, $J_n\in\mathcal{U}_n$.
This completes the proof of the claim.

As $L$ is not separable \cite{Stoyanova} Lemma 3.33 shows that $L^+\times L^-$ is not weakly Lindel\"of. $\Box$

Our third example is in a different direction: Spaces that are weakly Rothberger in finite powers need not have a weakly Menger product. Since the topological sum of two Rothberger spaces is Rothberger, and their product is a closed subspace of the square of their topological sum, the negation of Souslin's Hypothesis implies that there is a Rothberger space whose square is not weakly Lindel\"of. But if two Rothberger spaces are weakly Rothberger in their finite powers, must their product be weakly Rothberger?

\begin{theorem}[{\sf CH}]\label{weakmengerproduct} There are weakly Rothberger spaces $X$ and $Y$ such that
\begin{enumerate}
  \item{Each (finite or infinite) power of $X$ and of $Y$ is weakly Rothberger, and}
  \item{$X\times Y$ is not weakly Menger.}
\end{enumerate}
\end{theorem}
The proof of Theorem \ref{weakmengerproduct} is developed through a few propositions. Recall that for topological space $(X,{\mathcal T})$, $\textsf{PR}(X)$ denotes the collection of nonempty finite subsets
of $X$. For $S\in \textsf{PR}(X)$ and an open set $V\subseteq X$, $[S,V]$ denotes $\{T\in\textsf{PR}(X): S\subseteq T \subseteq V\}$.
The collection of subsets of the form $[S,V]$ of \textsf{PR}(X) is a basis for a topology, denoted $\textsf{PR}({\mathcal T})$ on $\textsf{PR}(X)$. Then $(\textsf{PR}(X), \textsf{PR}({\mathcal T}))$
is the Pixley–Roy space of $X$
If $X$ has a countable base, then $\textsf{PR}(X)\backslash\{\emptyset\}$
is a union of countably many sets, each with the finite intersection property; this implies
that $\textsf{PR}(X)$ has countable cellularity. But countable cellularity is equivalent to: each
element of $\mathcal D$ has a countable subset which is in $\mathcal D$.

\begin{proposition}\label{sumtoproductPR} Let $X$ and $Y$ spaces. Then $\textsf{PR}(X)\times \textsf{PR}(Y)$ is homeomorphic to $\textsf{PR}(X\oplus Y)$.
\end{proposition}

\Proof The function $\Phi:\textsf{PR}(X)\times \textsf{PR}(Y)\rightarrow \textsf{PR}(X\oplus Y)$ defined by $\Phi ((F,G))=F\cup G$
is one-to-one and onto, continuous and open, and thus a homeomorphism. $\Box$

At this point it is convenient to introduce another family of open covers: A family $\mathcal{F}$ of subsets of an infinite set $S$ is said to be an $\omega$-\emph{cover}{\footnote{Note that we are deviating from standard usage of the term $\omega$-cover: We do not require that the cover be an \underline{open} cover of a space.}} of $S$ if $S$ is not a member of $\mathcal{F}$, yet for each finite subset $F$ of $S$ there is a member $U$ of $\mathcal{F}$ such that $F\subseteq U$. Let $X$ be a topological space.
\[
  \Omega =\{\mathcal{U}\in\open: \mathcal{U}\mbox{ is an $\omega$-cover of }X\}.
\]

\begin{proposition} {\rm \cite{Sch1} \label{X+Ynotmenger}({\sf CH})} There are separable metric spaces $X$ and $Y$ each with the property $\sone({\Omega},{\Omega})$, but their topological sum $Z = X \oplus Y $ does not have $\sfin(\open,\open)$.
\end{proposition}

\Proof In \cite{Sch1} {\sf CH} is used to construct subsets $X$ and $Y$ of $^{\omega}{\Bbb Z}$ such that each has the property $\sone(\Omega, \Omega)$, but $(X \cup Y) \oplus (X\cup Y ) = ^{\omega}{\Bbb Z}$.

As was noted in Theorem 3.9 of \cite{COC2}, a topological space has the Menger property in all finite powers if, and ony if, it has the property $\sfin(\Omega,\Omega)$. Since $\sfin(\open,\open)$ is preserved by continuous images, closed subsets, and countable unions, since $\sfin({\Omega},{\Omega})$ is equivalent to $\sfin(\open,\open)$ in all finite powers, and since $\sfin({\Omega},{\Omega})$ is preserved by closed subsets, finite powers, and continuous images, and since $^{\omega}{\Bbb Z}$ (which is homeomorphic to the set of irrational numbers) does not have $\sfin(\open,\open)$, it follows that $(X \cup Y )^2$ does not have $\sfin(\open,\open)$. Thus neither $X\cup Y$ nor $X \times Y$ has $\sfin(\open,\open)$. It follows that $X\oplus Y$ does not have the Menger property. $\Box$

\begin{theorem}[Daniels, Theorem 5B]\label{DanielsTh5B} If $Z$ is a metrizable space with the property $\sone(\Omega,\Omega)$, then ${\sf PR}(Z)$ is weakly Rothberger in each (finite or infinite) power.
\end{theorem}
Thus Theorem \ref{DanielsTh5B} implies that each of ${\sf PR}(X)$ and ${\sf PR}(Y)$ is weakly Rothberger in all powers.
Next apply the following theorem to $X\oplus Y$:
\begin{theorem}[Daniels, Theorem 2A]\label{DanielsTh2A} If for a space $Z$, ${\sf PR}(Z)$ is weakly Menger, then each finite power of $Z$ is Menger.
\end{theorem}
Thus Theorem \ref{DanielsTh2A} states that if ${\sf PR}(X)$ has property $\sfin(\open,\dense)$, then $X$ has property $\sfin(\Omega,\Omega)$: For metrizable spaces the converse is implied by Theorem \ref{DanielsTh5B}.

It follows that ${\sf PR}(X\oplus Y)$ is not weakly Menger. But then Proposition \ref{sumtoproductPR} implies that ${\sf PR}(X)\times {\sf PR}(Y)$ is not weakly Menger. This completes the proof of Theorem \ref{weakmengerproduct}. $\Box$

As far as {\sf ZFC} results are concerned, there is also the following result by Todor\v{c}evic. First we introduce another family of open covers:
A family $\mathcal{F}$ of subsets of an infinite set $S$ is said to be a $\gamma$-\emph{cover} of $S$ if for each $x\in S$ the set $\{F\in {\mathcal F}:x\notin F\}$ is finite and $\mathcal F$ is infinite. Then for a topological space $X$ we define
\[
  \Gamma =\{\mathcal{U}\in\open: \mathcal{U}\mbox{ is an $\gamma$-cover of }X\}.
\]
Following Gerlits and Nagy, call a space which satisfies $\sone(\Omega,\Gamma)$ is a \emph{$\gamma$-space} \cite{GN}.

\begin{theorem}[Todor\v{c}evic, \cite{T}, Theorem 8]\label{Clines} There are ${\sf T}_3$ $\gamma$-spaces $X$ and $Y$ such that $X\times Y$ is not Lindel\"of.
\end{theorem}

And finally for this section:
\begin{theorem}\label{rothbsquare} It is consistent, relative to the consistency of {\sf ZFC}, that there is a ${\sf T}_3$-Rothberger space whose square is not Lindel\"of.
\end{theorem}
\Proof
Let ${\mathbb S}$ denote the \emph{Sorgenfrey line}, the topological space obtained from refining the standard topology on the real line by also declaring each interval of the form $\lbrack a,\, b)$ open. It is well known that $\sorgenfrey$ is a ${\sf T}_3$ Lindel\"of space while $\sorgenfrey\times\sorgenfrey$ is not Lindel\"of, but still ${\sf T}_3$. By Theorem 11 of \cite{MSFT}, if $\kappa$ is an uncountable cardinal and if $\cohen(\kappa)$ denotes the Cohen forcing notion for adding $\kappa$ Cohen reals, then in the generic extension the ground model Sorgenfrey line is a Rothberger space. But proper forcing preserves not being Lindel\"of, and ${\Bbb S}$ in the generic extension by $\cohen(\kappa)$, the square of the ground model copy of $\sorgenfrey$ is not Lindel\"of. Since ${\sf T}_3$ is preserved, it follows that the square of a (almost) Rothberger space need not be (almost) Lindel\"of.
$\Box$

As an aside to the proof of Theorem \ref{rothbsquare}: In Lemma 17 of \cite{BaSc1} it was shown that ${\mathbb S}$ does not have the property $\sfin(\open,\open)$, and since ${\mathbb S}$ is ${\sf T}_3$, this means that ${\mathbb S}$ is not almost Menger. As noted in the proof of Theorem \ref{rothbsquare}, in the generic extension by uncountably many Cohen reals, the ground model version of $\sorgenfrey$ is Rothberger and thus Menger. Thus, proper forcing does not preserve being not Menger.

Also note that the ground model Sorgenfrey line remains a separable space in the generic extension, and thus a weakly Rothberger space in a strong sense: TWO has a winning strategy in the game $\gone^{\omega}(\open,\dense)$.

\section{Productively Lindel\"of spaces}    

For a topological property {\sf Q} that is inherited by closed sets we shall say that a space is \emph{productively }{\sf Q} if for any space $Y$ which has property {\sf Q}, also $X\times Y$ has property {\sf Q}. Note that if a space $X$ is productively ${\sf Q}$, then each finite power of $X$ is also productively ${\sf Q}$.

Much work has been done on characterizing the members of the class of productively Lindel\"of spaces.
We expand the basic problem of characterizing the class of productively Lindel\"of spaces to also characterizing the classes of spaces that are productive for covering properties that are relatives of the Lindel\"of property. Some basic problems emerge:
Assume that we have topological properties ${\sf Q}$ and ${\sf R}$ where ${\sf Q}$ implies ${\sf R}$. When is it the case that:
\begin{enumerate}
  \item{If a space is productively ${\sf R}$, then it is productively ${\sf Q}$?}
  \item{If a space is productively ${\sf Q}$, then it is productively ${\sf R}$?}
  \item{If the product of $X$ with every space of property ${\sf Q}$ has property ${\sf R}$, then is $X$ productively ${\sf R}$?}
\end{enumerate}

Compact spaces are productively Lindel\"of but need not be Rothberger spaces, and thus need not be productively Rothberger. Thus, spaces that are productive for one class of Lindel\"of spaces need not be productive for another.

\begin{center}{\bf Alster's Theorems and a property of Alster.}\end{center}

In \cite{Alster1} Alster proves the following interesting theorem:
\begin{theorem}[Alster]\label{alsterth} Consider the following statements about a space $X$:
\begin{enumerate}
  \item{$X$ is an Alster space.}
  \item{$X$ is productively Lindel\"of.}
\end{enumerate}
Then (1) implies (2). If $X$ is a space of weight at most $\aleph_1$, and if {\sf CH} holds, then also (2) implies (1).
\end{theorem}

Alster \cite{Alster1}, Lemma 1, also proved that the Alster spaces give an affirmative answer to Problem \ref{prodLind}.
\begin{theorem}[Alster]\label{countable powers}
   If $X$ is an Alster space, then $X^{\aleph_0}$ is a Lindel\"of space.
\end{theorem}

The following classes of covers are central to Alster's analysis of the productively Lindel\"of spaces:

\begin{quote}
$\mathcal{G}_{K}$: The family consisting of sets $\mathcal{U}$ where $X$ is not in $\mathcal{U}$, each element of $\mathcal{U}$ is a ${\sf G}_{\delta}$ set, and for each compact set $C\subset X$ there is a $U\in\mathcal{U}$ such that $C\subseteq U$.
\end{quote}

\begin{quote}
${\mathcal{G}}$: The family of all covers $\mathcal{U}$ of the space $X$ for which each element of $\mathcal{U}$ is a ${\sf G}_{\delta}$-set.
\end{quote}

In \cite{BKR} Theorem 4.5 it is proved that the product of finitely many Alster spaces is again an Alster space. On account of this fact it is useful to also consider the following class of covers of spaces:
\begin{quote}
${\mathcal{G}}_{\Omega}$: This is the set of covers $\mathcal{U}\in\mathcal{G}$ for which $X$ is not in $\mathcal{U}$, but for each finite set $F\subset X$ there is a $U\in\mathcal{U}$ such that $F\subseteq U$.
\end{quote}
Observe that $\mathcal{G}_K$ is a subset of $\mathcal{G}_{\Omega}$.

The connection of Alster's condition to selection principles will now be determined through a sequence of Lemmas, culminating in Theorem \ref{alsterselchar}:
\begin{lemma}\label{AlsterS1} For a topological space $X$ the following are equivalent:
\begin{enumerate}
  \item{$X$ is an Alster space.}
  \item{$X$ satisfies the selection principle $\sone(\mathcal{G}_K,\mathcal{G})$.}
  \item{$X$ satisfies the selection principle $\sone(\mathcal{G}_K,\mathcal{G}_{\Omega})$.}
\end{enumerate}
\end{lemma}
\Proof
(1)$\Rightarrow$(2): Suppose that $X$ is an Alster space and let $(\mathcal{U}_n:n\in\naturals)$ be a sequence of elements of $\mathcal{G}_K$. Define
\[
  \mathcal{U} = \{\bigcap_{n\in\naturals}U_n:\, (\forall n)(U_n\in\mathcal{U}_n)\}.
\]
Then $\mathcal{U}$ is a member of $\mathcal{G}_K$. Since $X$ is Alster, choose a countable subset $(V_n:n\in\naturals)$ of $\mathcal{U}$ which is a cover of $X$. For each $n$ write
\[
  V_n = \bigcap_{k\in\naturals} U^n_k
\]
where for each $k$, $U^n_k$ is an element of $\mathcal{U}_n$. Finally for each $n$ set $W_n = U^n_n$, an element of $\mathcal{U}_n$. Then $\{W_n:n\in\naturals\}$ is a cover of $X$ and thus a member of $\mathcal{G}$.

{\flushleft{Proof of (2)$\Rightarrow$(1):}} Let a member $\mathcal{U}$ of ${\mathcal G}_K$ be given. For each $n$, set $\mathcal{U}_n = \mathcal{U}$. Then apply $\sone(\mathcal{G}_K,\mathcal{G})$ to this sequence and select for each $n$ a $U_n\in\mathcal{U}_n$ such that $\{U_n:n\in\naturals\}$ is a cover of $X$.

{\flushleft{Proof of (2)$\Rightarrow$(3):}}

{\bf Claim 1:} $X$ satisfies the selection principle $\sone(\mathcal{G}_K,\mathcal{G})$ if and only if $X^2$ satisfies the selection principle $\sone(\mathcal{G}_K,\mathcal{G})$ (hence every finite power of $X$ satisfies $\sone(\mathcal{G}_K,\mathcal{G})$).

\Proof Follows since the finite power of Alster space is an Alster space.

{\bf Claim 2:} Each finite power of $X$ satisfies the selection principle $\sone(\mathcal{G}_K,\mathcal{G})$ if and only if $X$ satisfies the selection principle $\sone(\mathcal{G}_K,\mathcal{G}_{\Omega})$.

\Proof Let $({\mathcal U}_n:n\in\naturals)$ be a sequence of elements of ${\mathcal G}_K$. Then $({\mathcal W}_n:n\in\naturals)$, where ${\mathcal W}_n=\{(U_n)^\kappa:U\in {\mathcal U}_n\}$, is a sequence of elements of ${\mathcal G}_K$ of $X^\kappa$.
Then we can select by hypothesis $(U_n)^\kappa\in {\mathcal W}_n$ such that $\{(U_n)^\kappa:n\in \naturals\}$ is a cover of $X^\kappa$.
Consider now a finite subset $F=\{x_1,\cdots,x_\kappa\}$ of $X$. We consider $F$ like a point of $X^\kappa$, say $z=(x_1,\cdots,x_\kappa)$. So there is an element $(U_n)^\kappa$ such that $z\in (U_n)^\kappa$. Then there is an element $U_n\in {\mathcal U}_n$ such that $F\subset U_n$. It follows that $\{U_n:n\in\naturals\}$ witnesses that $X$ satisfies $\sone({\mathcal G}_K,{\mathcal G}_\Omega)$.
The converse is obvious. This complete the proof.

{\flushleft{Proof of (3)$\Rightarrow$(2):}} Obvious.
$\Box$

A space $X$ is has the \emph{Hurewicz property} if for each sequence $({\mathcal U}_n : n\in\naturals)$ of open covers of $X$ there is a sequence $({\mathcal V}_n : n\in\naturals)$ such that for each $n$, ${\mathcal V}_n$ is a finite subset of ${\mathcal U}_n$ and each $x\in X$ belongs to $\bigcup{\mathcal V}_n$ for all but finitely many $n$. The Alster property implies the following strengthening of the Hurewicz property:

\begin{lemma}
If $X$ is a space that has property $\sone(\mathcal{G}_K,\mathcal{G})$, then there is for each sequence $(\mathcal{U}_n:n\in\naturals)$ of elements of $\mathcal{G}_K$ a sequence $(\mathcal{V}_n:n\in\naturals)$ such that for each $n$ we have $\mathcal{V}_n\subseteq \mathcal{U}_n$, $\vert\mathcal{V}_n\vert\le n$, and for each $x\in X$, for all but finitely many $n$, $x\in\bigcup \mathcal{V}_n$.
\end{lemma}
\Proof Let a sequence $(\mathcal{U}_n:n\in\naturals)$ of elements of $\mathcal{G}_K$ be given. For each $n$ define $\mathcal{W}_n$ to be the set $\{\bigcap_{n\in\naturals} U_n:\, (\forall n)(U_n\in\mathcal{U}_n)\}$. Applying $\sone(\mathcal{G}_K,\mathcal{G})$ to the sequence $(\mathcal{W}_n:n\in\naturals)$ we find for each $n$ a $W_n\in\mathcal{W}_n$ such that for each $x\in X$ there is an $n$ with $x\in W_n$.

For each $n$ write $W_n = \bigcap U^n_k$ where for each $n$ and $k$ we have $U^n_k\in\mathcal{U}_k$. Then, for each $k$, set
\[
  \mathcal{V}_k = \{U^1_k,\,\cdots,\,U^k_k\},
\]
a finite subset of $\mathcal{U}_k$. Note that if $x\in X$ is an element of $W_n$, then for each $k\ge n$ we have $x\in \cup\mathcal{V}_k$.
$\Box$

For each $n$, let $T_n$ be the $n$-th triangular number\footnote{$T_n=1+2+\cdots+n$.}. Using the technique in the proof of the previous Lemma, we find
\begin{lemma}\label{grouping}
If $X$ is a space that has property $\sone(\mathcal{G}_K,\mathcal{G})$, then there is for each sequence $(\mathcal{U}_n:n\in\naturals)$ of elements of $\mathcal{G}_K$ a sequence $(U_n:n\in\naturals)$ such that for each $n$ we have $U_n\in \mathcal{U}_n$, and for each $x\in X$, for all but finitely many $n$,
\[
  x\in \bigcup_{T_n< j\le T_{n+1}}U_j.
\]
\end{lemma}
\Proof Let a sequence $(\mathcal{U}_n:n\in\naturals)$ of elements of $\mathcal{G}_K$ be given. For each $n$ define
\[
  \mathcal{V}_n = \{\bigcap_{T_n<i\le T_{n+1}} U_i:\, (\forall i)(T_n<i\le T_{n+1})(U_i\in\mathcal{U}_i)\}.
\]
Now apply the conclusion of the previous lemma to the sequence $(\mathcal{V}_n:n\in\naturals)$ to find for each $n$ a set $\mathcal{W}_n\subseteq\mathcal{V}_n$ of cardinality $n$ such that for each $x$, for all but finitely many $n$, $x\in\bigcup\mathcal{W}_n$. Each element of $\mathcal{W}_n$ is of the form $U^j_{T_n+1}\cap\cdots\cap U^j_{T_{n+1}}$ where $1\le j\le n$ and each $U^j_i$ is an element of $\mathcal{U}_i$.

Now for each $m$ choose $V_m\in\mathcal{U}_m$ as follows: Find the largest $n$ with $T_n<m$ and then identify $j$ with $1\le j\le n+1$ with $m=T_n+j$, and put $V_m = U^j_{T_n+j}$. $\Box$

\begin{quote}
${\mathcal{G}^{gp}}$: This is the set of covers $\mathcal{U}\in\mathcal{G}$ for which there is a (disjoint) partition $\mathcal{U} = \bigcup_{n\in\naturals}\mathcal{U}_n$ such that each $\mathcal{U}_n$ is finite, and for each $x\in X$ for all but finitely many $n$, $x$ is in $\bigcup\mathcal{U}_n$.
\end{quote}

\begin{theorem}\label{alsterselchar} For a topological space $X$ the following are equivalent:
\begin{enumerate}
  \item{$X$ is an Alster space.}
  \item{$X$ satisfies the selection principle $\sone(\mathcal{G}_K,\mathcal{G})$.}
  \item{$X$ satisfies the selection principle $\sone(\mathcal{G}_K,\mathcal{G}_{\Omega})$.}
  \item{$X$ satisfies the selection principle $\sone(\mathcal{G}_K,\mathcal{G}^{gp})$.}
  \item{$X$ satisfies the selection principle $\sone(\mathcal{G}_K,\mathcal{G}_{\Omega}^{gp})$.}
  \item{$X$ satisfies the selection principle $\sfin(\mathcal{G}_K,\mathcal{G})$.}
  \item{$X$ satisfies the selection principle $\sfin(\mathcal{G}_K,\mathcal{G}^{gp})$.}
  \item{$X$ satisfies the selection principle $\sfin(\mathcal{G}_K,\mathcal{G}_{\Omega}^{gp})$.}
\end{enumerate}
\end{theorem}

\begin{corollary}
If $X$ is an Alster space then $X$ has a Hurewicz space in all finite powers.
\end{corollary}

\begin{center}{\bf Strengthening Alster's property: $\sone(\mathcal{G}_K,\mathcal{G}_{\Gamma})$.}\end{center}

Of several standard ways in which to strengthen the Alster property, consider the following one. Define
\begin{quote}
${\mathcal{G}}_{\Gamma}$: This is the set of covers $\mathcal{U}\in\mathcal{G}$ which are infinite, and each infinite subset of $\mathcal{U}$ is a cover of $X$.
\end{quote}
Then $\sone(\mathcal{G}_K,\mathcal{G}_{\Gamma})$ is a formal strengthening of the Alster property. It is not clear if this formal strengthening is in fact a real strengthening. Several known examples of productively Lindel\"of spaces are so because they have this stronger version of Alster's property. We review some of these:

{\flushleft{{\bf A:}}} An easy consequence of one of Alster's results gives:
\begin{theorem}[{\sf CH}]\label{AlsterTh3} Assume that $X$ is a space of weight at most $\aleph_1$, and that each compact subset of $X$ is a ${\sf G}_{\delta}$ set. If $X$ has the property $\sone(\mathcal{G}_K,\mathcal{G})$, then $X$ has the property $\sone(\mathcal{G}_K,\mathcal{G}_{\Gamma})$.
\end{theorem}
\Proof It follows immediately from Alster's Theorem, Theorem \ref{alsterth}, that if $X$ is a space of weight at most $\aleph_1$ in which each compact set is ${\sf G}_{\delta}$  and if $X$ is productively Lindel\"of, then {\sf CH} implies that $X$ is $\sigma$-compact (as the set of compact subsets of $X$ is a cover of $X$ by ${\sf G}_{\delta}$ sets of the required kind). Thus, under {\sf CH}, every productively Lindel\"of space of weight at most $\aleph_1$ for which each compact subspace is a ${\sf G}_{\delta}$ set has the property $\sone(\mathcal{G}_K,\mathcal{G}_{\Gamma})$.
$\Box$

{\flushleft{{\bf B:}}} A topological space is said to be a {\sf P}-\emph{space} if each intersection of countably many open sets is open. Galvin (see \cite{GN}) pointed out that any Lindel\"of {\sf P}-space is a $\gamma$-space. Thus, as the topology of a {\sf P}-space is the ${\sf G}_{\delta}$ topology, we find:
\begin{lemma}[Galvin]\label{GalvinLemma} Each Lindel\"of {\sf P}-space has the property $\sone(\mathcal{G}_K,\mathcal{G}_{\Gamma})$.
\end{lemma}
In \cite{KF} Proposition 2.1 the authors show  that the product of a Lindel\"of {\sf P}-space with any Lindel\"of space is a Lindel\"of space. This result now follows from (1) $\Rightarrow$ (2) of Theorem \ref{alsterth}, and Lemma \ref{GalvinLemma}. A classical theorem of Noble \cite{Noble} states that the product of countably many Lindel\"of {\sf P}-spaces is still Lindel\"of. This result now follows from Theorem \ref{countable powers} and Lemma \ref{GalvinLemma}. Additionally it is also known that:
\begin{theorem}[Misra, \cite{Misra}]\label{PRothb}
A Lindel\"of {\sf P}-space is productively Lindel\"of {\sf P}.
\end{theorem}
\bigskip

{\flushleft{{\bf C:}}} A space is \emph{scattered} if each nonempty subspace has an isolated point.
\begin{theorem}[Gewand \cite{Gewand}, Theorem 2.2]\label{GewandTh}
If $X$ is a scattered Lindel\"of space, then in the ${\sf G}_{\delta}$ topology $X$ is a Lindel\"of {\sf P}-space.
\end{theorem}
\begin{corollary}\label{scattered2} A scattered Lindel\"of space satisfies $\sone(\mathcal{G}_K,\mathcal{G}_{\Gamma})$.
\end{corollary}
It follows that Lindel\"of scattered spaces are productively Lindel\"of, and that the countable power of such a space is Lindel\"of.

{\flushleft{{\bf D:}}} A space is $\sigma$-compact if it is the union of countably many compact subsets.
\begin{theorem}\label{sigmacompact} Each $\sigma$-compact space has the property $\sone(\mathcal{G}_K,\mathcal{G}_{\Gamma})$.
\end{theorem}
\Proof Let $X$ be $\sigma$-compact and write $X = \bigcup_{n\in\naturals} K_n$ where for each $n$ we have $K_n\subseteq K_{n+1}$ and $K_n$ is compact. Let a sequence $(\mathcal{U}_n:n\in\naturals)$ of elements of $\mathcal{G}_K$ be given. For each $n$, choose $U_n\in\mathcal{U}_n$ with $K_n\subseteq U_n$. Then each $U_n$ is a ${\sf G}_{\delta}$ subset of $X$, and for each $x\in X$, for all but finitely many $n$, $x\in U_n$. $\Box$

This implies the well-known fact that $\sigma$-compact spaces are productively Lindel\"of. Via Theorem \ref{countable powers} we also find that the countable power of a $\sigma$-compact space is Lindel\"of.

We now explore the productive properties of spaces with this stronger version of the Alster property.
\begin{theorem}\label{AlsterHurewicz} If $X$ is a topological space with property $\sone(\mathcal{G}_K,\mathcal{G}_{\Gamma})$, then $X$ is:
\begin{enumerate}
  \item{productively Lindel\"of.}
  \item{productively Menger.}
  \item{productively Hurewicz.}
\end{enumerate}
\end{theorem}
\Proof We show that $X$ is productively Hurewicz. The proof that $X$ is productively Menger is similar.

Let $Y$ be a Hurewicz space, and let $(\mathcal{U}_n:n\in\naturals)$ be a sequence of open covers of $X\times Y$. We may assume that each $\mathcal{U}_n$ is closed under finite unions. For each compact subset $K$ of $X$, and for each $n$, find a ${\sf G}_{\delta}$ set $\phi_n(K)\supset K$ such that for each $y\in Y$ there is an open set $U\in\mathcal{U}_n$ with $\phi(K)\times\{y\}\subseteq U$.

This defines, for each $n$, a set $\mathcal{G}_n =\{\phi_n(K):K\subset X \mbox{ compact}\} \in\mathcal{G}_K$. Applying $\sone(\mathcal{G}_K,\mathcal{G}_{\Gamma})$ to $(\mathcal{G}_n:n\in\naturals)$, we find for each $n$ a compact set $K_n\subset X$ such that $\{\phi_n(K_n):n\in\naturals\}$ is a member of $\mathcal{G}_{\Gamma}$.

Now for each $n$ define
\[
  \mathcal{H}_n =\{V\subset Y:\, V\mbox{ open and there is a } U\in\mathcal{U}_n \mbox{ with }\phi_n(K_n)\times V\subseteq U\}.
\]

Apply the fact that $Y$ is a Hurewicz space to the sequence $(\mathcal{H}_n:n\in\naturals)$: For each $n$ choose a finite set $\mathcal{J}_n\subseteq \mathcal{H}_n$ such that for each $y\in Y$, for all but finitely many $n$, $y\in\cup\mathcal{J}_n$.

Finally, for each $n$ choose for each $H\in\mathcal{J}_n$ a $U_H\in\mathcal{U}_n$ with $\phi_n(K_n)\times H\subseteq U_H$, and put $\mathcal{V}_n = \{U_H:\, H\in\mathcal{J}_n\}$. Then we have for each $n$ that $\mathcal{V}_n$ is a finite subset of $\mathcal{U}_n$, and for each $(x,y)\in X\times Y$, for all but finitely many $n$, $(x,y)$ is a member of $\bigcup\mathcal{V}_n$.
$\Box$

There are many {\sf ZFC} examples of Lindel\"of {\sf P}-spaces. For example in \cite{MSRBGs} it is shown that for each uncountable cardinal number $\kappa$ there is a ${\sf T}_{3\frac{1}{2}}$ Lindel\"of {\sf P} -group of cardinality $\kappa$ on which TWO has a winning strategy in the game $\gone^{\omega}(\Omega,\Gamma)$. Since the topology is the ${\sf G}_{\delta}$ topology, this means that in these examples TWO has a winning strategy in the game $\gone(\mathcal{G}_K,\mathcal{G}_{\Gamma})$.

Now we need two more classes of open covers: A family $\mathcal{F}$ of subsets of an infinite set $S$ is said to be a \emph{large cover} of $S$ if for each $x\in S$ the set $\{F\in {\mathcal F}:x\in F\}$ is infinite.

A family $\mathcal F$ is a \emph{groupable} cover if, and only if, there is a partition ${\mathcal F} =\bigcup_{n<\infty} U_n$ where the $U_n$'s are finite and disjoint from each other, such that each point in the space belongs to all but finitely many of the sets $U_n$.

We associate the following symbols with classes of open covers that are large or groupable:
\[
  \Lambda :=\{\mathcal{U}\in\open: \mathcal{U}\mbox{ is an \emph{large}-cover of }X\}.
\]
\[
  \open^{gp}:=\{\mathcal{U}\in\open: \mathcal{U}\mbox{ is groupable in }X\}.
\]
A space $X$ is called a \emph{Gerlits-Nagy space} if it satisfies the selection principle $\sone(\Omega, \open^{gp})$ \cite{COC7}. Each $\gamma$-space is a Gerlits -Nagy space, and each Gerlits-Nagy space is a Rothberger space.

\begin{theorem}\label{AlsterRothb} If $X$ is a Rothberger space with property $\sone(\mathcal{G}_K,\mathcal{G}_{\Gamma})$, then $X$ is:
\begin{enumerate}
  \item{productively Rothberger.}
  \item{productively Gerlits-Nagy.}
\end{enumerate}
\end{theorem}
\Proof (of (1)) Let a Rothberger space $Y$ be given, and let $(\mathcal{U}_n:n\in\naturals)$ be a sequence of open covers of $X\times Y$. As $X$ is Rothberger, each of its compact subsets is Rothberger. Write $\naturals = \bigcup_{n\in\naturals}S_n$ where the $S_n$'s are infinite and pairwise disjoint.

Fix $n$, and for each compact subset $C$ of $X$, and for each $y\in Y$, choose a finite sequence $U_{i_1},\cdots,U_{i_k}$ where
\begin{enumerate}
  \item{$i_1<\cdots< i_k$ are elements of $S_n$ and}
  \item{$U_{i_j}$ is an element of $\mathcal{U}_{i_j}$, $1\le j\le k$, and}
  \item{$C\times \{y\}\subseteq U_{i_1}\cup\cdots\cup U_{i_k}$.}
\end{enumerate}
This is possible as compact subsets of $X$ are Rothberger spaces.

Thus, for fixed $n$, for each compact subset $C$ of $X$ we find that $C\times Y$ is Rothberger, thus Lindel\"of, and we can find a ${\sf G}_{\delta}$ subset $\phi_n(C)$ of $X$ such that $C\subseteq \phi_n(C)$, and for each $y\in Y$ there is a sequence $i_1<\cdots< i_k$ of elements of $S_n$ such that $\phi_n(C)\times\{y\}$ is a subset of $U_{i_1}\cup\cdots\cup U_{i_k}$ as above. But then $\mathcal{G}_n = \{\phi_n(C):C\subset X\mbox{ compact}\}$ is a member of $\mathcal{G}_K$ for $X$.

Apply $\sone(\mathcal{G}_K,\mathcal{G}_{\Gamma})$ to the sequence $(\mathcal{G}_n:n\in\naturals)$ and select for each $n$ a $G_n\in\mathcal{G}_n$ such that for each $x\in X$, for all but finitely many $n$, $x\in G_n$.

For fixed $n$, choose for each $y\in Y$ an open set $U^n_y\subset Y$ such that $y\in U^n_y$ and there is a sequence $i_1<\cdots< i_k$ of elements of $S_n$ such that $G_n\times U^n_y$ is a subset of $U_{i_1}\cup\cdots\cup U_{i_k}$ as above. Then $\mathcal{H}_n = \{U^n_y:y\in Y\}$ is an open cover of $Y$.

Next apply the fact that $Y$ is Rothberger to the sequence $(\mathcal{H}_n:n\in\naturals)$, and choose for each $n$ an $H_n\in\mathcal{H}_n$ such that for each $y\in Y$ there are infinitely many $n$ with $y\in H_n$.

Then for each $n$ choose $i^n_1<\cdots< i^n_{k_n}\in S_n$ such that $G_n\times H_n\subseteq U_{i^n_1}\cup\cdots\cup U_{i^n_{k_n}}$.

It follows that there is a sequence of $U_n\in\mathcal{U}_n$ such that $\{U_n:n\in\naturals\}$ is an open cover of $X\times Y$.

The poof of (2) is similar, and left to the reader. $\Box$

It follows that Lindel\"of {\sf P}-spaces, Lindel\"of scattered spaces, as well as $\sigma$-compact Rothberger spaces are not only productively Lindel\"of, but also productively Menger, productively Hurewicz, productively Rothberger and productively Gerlits-Nagy.

Regarding the hypothesis of Theorem \ref{AlsterRothb} that $X$ should be Rothberger: Note that the hypothesis $\sone(\mathcal{G}_K,\mathcal{G})$ implies $\sone(\open_K,\open)$, which implies $\sfin(\open,\open)$. In Example 2 of \cite{BScountabledim} it was shown that {\sf CH} implies the existence of a subspace $X$ of $\reals^{\naturals}$ which is not countable dimensional and yet satisfies $\sone(\open_K,\Gamma)$. Such an $X$ cannot be a Rothberger space, since metrizable Rothberger spaces are zero dimensional. One may consider weakening the hypothesis that $X$ is Rothberger to the hypothesis that each compact subset of $X$ is Rothberger. However, this is not a weakening of the hypotheses, since:
\begin{lemma}\label{Rothb} If $X$ is a space such that each compact subspace is Rothberger, and $\sone(\open_K,\open)$ holds of $X$ then $X$ is a Rothberger space.
\end{lemma}

Similarly, one can show:
\begin{lemma}\label{Rothb1} If $X$ is a space such that each compact subspace is Rothberger, and $\sone(\open_K,\open^{gp})$ holds of $X$ then $X$ is a Gerlits-Nagy space.
\end{lemma}

We leave the proof of these lemmas to the reader. We do not know if a similar statement is true about $\gamma$-spaces:
\begin{problem} Is it true that if each compact subset of $X$ is a Rothberger space and $X$ has the property $\sone(\open_K,\Gamma)$, then $X$ has the property $\sone(\Omega,\Gamma)$?
\end{problem}
What we can prove is:
\begin{lemma}\label{Rothb2} If $X$ is a space such that each compact subspace is finite, and $\sone(\open_K,\Gamma)$ holds of $X$ then $X$ is a $\sone(\Omega,\Gamma)$-space.
\end{lemma}

We now record the few results we have on productively $\gamma$ spaces.
\begin{lemma}{\rm\cite{COC2}}\label{omegacoverrefine} Every open $\omega$-cover of $X\times Y$ is refined by one whose elements are of the form $U\times V$ where $U\subseteq X$ and $V\subseteq Y$ are open.
\end{lemma}

\begin{theorem}\label{gammaproductive}
If $X$ is a space satisfying $\sone(\mathcal{G}_K,\mathcal{G}_{\Gamma})$ and if each compact subset of $X$ is finite, then $X$ is $\sone(\Omega,\Gamma)$-productive.
\end{theorem}
\Proof We know that $X$ is productively Lindel\"of. Let a $\gamma$-space $Y$ be given, and let $(\mathcal{U}_n:n\in\naturals)$ be a sequence of open $\omega$-covers of $X\times Y$.

{\flushleft{\bf Claim}} $X\times Y$ is Lindel\"of in each finite power.\\
For $\gamma$-spaces it is well known that each finite power of a $\gamma$-space is a $\gamma$-space. Similarly for {\sf P}-spaces, each finite power of a {\sf P}-space is a {\sf P}-space. Moreover, $(X\times Y)^n$ is homeomorphic to $X^n\times Y^n$.
Note that since $X\times Y$ is a Lindel\"of space in each finite power, we may assume, by the Proposition on p. 156 of \cite{GN}, that each $\mathcal{U}_n$ is a countable set. For each $n$ we may also assume by Lemma \ref{omegacoverrefine} that  the elements of $\mathcal{U}_n$ are of the form $U\times V$ where $U$ is open in $X$ and $V$ is open in $Y$. Thus, each $\mathcal{U}_n$ is of the form $\{U^n_k\times V^n_k:k\in\naturals\}$.

Fix $n$ and fix a finite subset $F$ of $X$. For each finite set $G\subset Y$ choose a $U^n_{n(F,G)}\times V^n_{n(F,G)}\in\mathcal{U}_n$ containing $F\times G$. The intersection $W^n(F)=\bigcap\{U_{n(F,G)}: G\subset Y\mbox{ finite}\}$ is a countable intersection as $\mathcal{U}_n$ is countable, and thus is an element of $\mathcal{G}_K$ for $X$. Likewise, $W(F) = \cap_{n\in\naturals}W^n(F)$ is and element of $\mathcal{G}_K$ for $X$, and for each finite $G\subset Y$ and each $n$, $W(F)\times V^n_{n(F,G)}$ contains $F\times G$.

Note that $\{W(F):F\subset X \mbox{ finite}\}$ is an $\omega$-cover of $X$. Since $X$ is a  Lindel\"{o}f space satisfying $\sone(\mathcal{G}_K,\mathcal{G}_{\Gamma})$, choose a countable set $\{F_n:n\in\naturals\}$ of finite subsets of $X$ with $\{W(F_n):n\in\naturals\}$ an element of $\mathcal{G}_{\Gamma}$.

Fix $n$: Then $\{W(F_n)\times V^n_k: k\in\naturals \mbox{ and } W(F_n)\subseteq U^n_k\}$ is an $\omega$-cover of $W(F_n)\times Y$, and also ${\mathcal Z}_n = \{V^n_k:k\in\naturals \mbox{ and }W(F_n)\subseteq U^n_k\}$ is an $\omega$-cover for $Y$.

Since $Y$ is a $\gamma$-set, choose for each $n$ a $k_n$ such that $\{V^n_{k_n}:n\in\naturals\}$ is a $\gamma$-cover of $Y$.
Then $\{W(F_n)\times V^n_{k_n}:n\in\naturals\}$ is a $\gamma$-cover for $X\times Y$, and so $\{U^n_{k_n}\times V^n_{k_n}: n\in\naturals\}$ is a $\gamma$-cover of $X\times Y$. $\Box$

Since in a ${\sf T}_2$ Lindel\"of {\sf P}-space compact subsets are finite, it follows that ${\sf T}_2$ Lindel\"of {\sf P}-spaces as well as ${\sf T}_2$ Lindel\"of scattered spaces are productively $\gamma$-spaces. As ${\sf T}_2$ compact Rothberger spaces are scattered Lindel\"of spaces, these are productively $\gamma$-spaces.

\begin{corollary} $\sigma$-compact Rothberger spaces are productively $\gamma$.
\end{corollary}
\Proof Since $X$ is $\sigma$-compact we write $X =\bigcup_{n\in\naturals}X_n$ where each $X_n$ is a compact Rothberger space. Since the union of finitely many compact Rothberger spaces is a compact Rothberger space, we may assume that for each $n$ we have $X_n\subseteq X_{n+1}$. A compact Rothberger space is a scattered Lindel\"of space. Then for $Y$ a $\gamma$-space, for each $n$  $X_n\times Y$ is a $\gamma$-space. This gives a union of an increasing sequence of $\gamma$-spaces, and so by Jordan's theorem (\cite{FJordan}, Corollary 14) is again a $\gamma$-space. $\Box$

By Corollary 15 of \cite{MSRBGs} there is for each infinite cardinal $\kappa$ a $\sigma$ -compact ${\sf T}_0$ topological group of cardinality $\kappa$ such that TWO has a winning strategy in the game $\gone(\Omega,\Gamma)$. By Corollary 17 of \cite{MSRBGs} there is for each infinite cardinal number $\kappa$ a ${\sf T}_0$ topological group $(G, *)$ of cardinality $\kappa$ which is a $\sigma$ -compact Rothberger space in all finite powers.

Alster also showed that
\begin{theorem}[Alster, \cite{Alster1}, Theorem 4]\label{Aslterth4} Assume {\sf CH}. If $X$ is of weight at most $\aleph_1$ and has property $\sone(\mathcal{G}_K,\mathcal{G})$, and if each compact subset of $X$ is at most countable, then in the ${\sf G}_{\delta}$ topology $X$ is Lindel\"of.
\end{theorem}
It follows that in addition to be productively Lindel\"of such spaces are also productively Menger-, Hurewicz-, Rothberger-, Gerlits-Nagy- and $\gamma$.

\section{Productivity of weak covering properties.}

Spaces that are productively weakly Lindel\"of do not currently have as well developed a theory. A number of ways of generalizing the notion of an Alster space or its strengthenings as considered in the previous section suggest themselves, but we have had limited success in exploiting these to identify classes of spaces that are for example productively weakly Lindel\"of spaces. In this section we report some of these results, and pose a number of questions whose answers may help identify criteria under which a space with a weak covering property is productively so.

The following two elementary properties of the notion of dense set will be used several times.
\begin{lemma}\label{densetransitive} Let $X$ be a topological space. If $Y\subset X$ is dense in $X$, and $D\subset Y$ is dense in $Y$, then $D$ is dense in $X$.
\end{lemma}

\begin{lemma}\label{productdense} If $D\subset X$ is dense in $X$ and $E\subset Y$ is dense in $Y$, then $D\times E$ is dense in $X\times Y$.
\end{lemma}

The power of these two lemmas lie in the following:
\begin{theorem}\label{denseupwabsolute} Let $X$ be a topological with dense subset $D$.
\begin{enumerate}
  \item{If $D$ is productively weakly Lindel\"of, so is $X$.}
  \item{If $D$ is productively weakly Menger, so is $X$.}
  \item{If $D$ is productively weakly Hurewicz, so is $X$.}
  \item{If $D$ is productively weakly Rothberger, so is $X$.}
  \item{If $D$ is productively weakly Gerlits-Nagy, so is $X$.}
\end{enumerate}
\end{theorem}
\Proof We show the argument for productively weakly Rothberger, leaving the rest to the reader.
Thus, let $Y$ be a weakly Rothberger space and let $D$ be productively weakly Rothberger. Let $(\mathcal{U}_n:n\in\naturals)$ be a sequence of open covers of $X\times Y$. Then the relativizations to $D\times Y$ is a sequence of open covers of $D\times Y$, since $D$ is dense in $X$. Applying the fact that $D\times Y$ is weakly Rothberger we find in each $\mathcal{U}_n$ an element $U_n$ such that $\bigcup_{n\in\naturals}U_n$ is dense in $D\times Y$, and so by Lemma \ref{densetransitive} is dense in $X\times Y$.
$\Box$

Also the following fact is useful:
\begin{theorem}\label{refinement} Let $(X,\tau)$ be a topological and let $\tau^{\prime}$ be a finer topology on $X$ (i.e., $\tau\subset  \tau^{\prime}$).
\begin{enumerate}
  \item{If $(X,\tau^{\prime})$ is productively weakly Lindel\"of, so is $(X,\tau)$.}
  \item{If $(X,\tau^{\prime})$ is productively weakly Menger, so is $(X,\tau)$.}
  \item{If $(X,\tau^{\prime})$ is productively weakly Hurewicz, so is $(X,\tau)$.}
  \item{If $(X,\tau^{\prime})$ is productively weakly Rothberger, so is $(X,\tau)$.}
  \item{If $(X,\tau^{\prime})$ is productively weakly Gerlits-Nagy, so is $(X,\tau)$.}
\end{enumerate}
\end{theorem}

In what follows we find conditions under which a space is productively weakly {\sf T}, where {\sf T} is one of the weak covering properties we are considering.

\begin{center}{\bf Weakly compact spaces.}
\end{center}

The space $(X,\tau)$ is \emph{weakly} compact if there is for each open cover of the space a finite subset with union dense in the space. Since the closure of a finite union of sets is the union of the finitely many sets individual closures, the weakly compact spaces coincide with the almost compact spaces, and in the context of ${\sf T}_2$-spaces, coincide with the {\sf H}-closed spaces.

Analogous to Tychonoff's theorem for compact spaces, one has:
\begin{theorem}[Scarborough and Stone \cite{SS} Theorem 2.4] The product of weakly compact spaces is weakly compact.
\end{theorem}

The argument in Proposition 1.9 of \cite{CS} shows
\begin{theorem}\label{weakcompact} If $X$ is a weakly compact space then $X$ is productively weakly Lindel\"of. \end{theorem}
This can be extended to the following:
\begin{theorem}
Let $X$ be a weakly compact space.
\begin{enumerate}
  \item{$X$ is productively weakly Menger.}
  \item{$X$ is productively weakly Hurewicz.}
\end{enumerate}
\end{theorem}

Since compact spaces are weakly compact, these results also imply that compact spaces are productively weakly Lindel\"of, productively weakly Menger, and productively weakly Hurewicz.

\begin{center}{\bf Weak versions of the Alster property}
\end{center}

Towards exploration of productiveness for the weaker versions of the covering properties we considered, we introduce the following family of subsets of a topological space:
\begin{quote}
${\mathcal G}_D$ denotes the collection of sets $\mathcal{U}$ where each element of $\mathcal{U}$ is a ${\sf G}_{\delta}$ set, and $\bigcup\mathcal{U}$ is dense in the space.
\end{quote}

A space is said to be \emph{weakly Alster} if each member of $\mathcal{G}_{K}$ has a countable subset which is a member of $\mathcal{G}_D$.

Using the proof technique of Theorem 4.5 of \cite{BKR} we find:
\begin{lemma}\label{weakalsterLindelof}
If $X$ weakly Alster and $Y$ is weakly Lindel\"of, then $X\times Y$ is weakly Lindel\"of.
\end{lemma}
\Proof
Let $\mathcal{U}$ be an open cover of $X\times Y$. We may assume that $\mathcal{U}$ is closed under finite unions. For each compact set $C\subset X$ and for each $y\in Y$ we find an open set $U(C,y)\in \mathcal{U}$ such that $C\times \{y\}\subseteq U(C,y)$.  For each $y$ we find open sets $V(C,y)\subset X$ and $W(C,y)\subset Y$ such that $C\times\{y\}\subset V(C,y)\times W(C,y)\subset U(C,y)$.

By Theorem \ref{weakcompact} $C\times Y$ is weakly Lindel\"of. Thus the cover $\{W(C,y):y\in Y\}$ of $Y$ contains a countable subset with union dense in $Y$, say $\{W(C,y^C_n)
:n\in\naturals\}$. Define $V(C) = \bigcap\{V(C,y^C_n):n\in\naturals\}$, a ${\sf G}_{\delta}$ subset of $X$ containing the compact set $C\subseteq X$. But then $\mathcal{V} = \{V(C):C\subset X \mbox{ compact}\}$ is a member of $\mathcal{G}_K$ for $X$. Since $X$ is weakly Alster we find a countable subset $\{V(C_m):m\in\naturals\}$ of $\mathcal{V}$ with union dense in $X$. But then
$\{W(C_m,y^{C_m}_n): m,n\in\naturals\} \subset \mathcal{U}$ is countable and its union is dense in $X\times Y$.
$\Box$

Since an Alster space is a weakly Alster space, we find
\begin{corollary}[{\sf CH}]\label{plpwl} Every productively Lindel\"of space of weight at most $\aleph_1$ is productively weakly Lindel\"of.
\end{corollary}
\Proof By Theorem \ref{alsterth}, productively Lindel\"of spaces of weight at most $\aleph_1$ are Alster spaces, and thus weakly Alster spaces. $\Box$

We don't know if the additional hypotheses are necessary in Corollary \ref{plpwl}:
\begin{problem}\label{lindprodweaklindprod} Is every productively Lindel\"of space productively weakly Lindel\"of?
\end{problem}

\begin{problem} Is every productively weakly Lindel\"of space a weakly Alster space?
\end{problem}

An argument similar to the one in Theorem \ref{weakalsterLindelof} shows
\begin{theorem}\label{weakalsterproducts}
If $X$ and $Y$ are weakly Alster, then $X\times Y$ is weakly Alster.
\end{theorem}

Additionally, analogous to Lemma \ref{AlsterS1} we find
\begin{lemma}\label{weakAlsterS1} For a topological space $X$ the following are equivalent:
\begin{enumerate}
  \item{$X$ is a weakly Alster space.}
  \item{$X$ satisfies the selection principle $\sone(\mathcal{G}_K,\mathcal{G}_{D})$.}
  \item{$X$ satisfies the selection principle $\sone(\mathcal{G}_K,\mathcal{G}_{D_{\Omega}})$.}
  \item{$X$ satisfies the selection principle $\sone(\mathcal{G}_K,\mathcal{G}_{D^{gp}})$.}
  \item{$X$ satisfies the selection principle $\sone(\mathcal{G}_K,\mathcal{G}_{D_{\Omega}^{gp}})$.}
  \item{$X$ satisfies the selection principle $\sfin(\mathcal{G}_K,\mathcal{G}_D)$.}
  \item{$X$ satisfies the selection principle $\sfin(\mathcal{G}_K,\mathcal{G}_{D^{gp}})$.}
  \item{$X$ satisfies the selection principle $\sfin(\mathcal{G}_K,\mathcal{G}_{D_{\Omega}^{gp}})$.}
\end{enumerate}
\end{lemma}

We also leave to the reader the proof of
\begin{lemma}\label{densegoodsubspaces}
If the ${\sf T}_2$-space $X$ has a dense subspace which is a weakly Alster space, then $X$ is a weakly Alster space.
\end{lemma}

\begin{center}{\bf A stronger version of weakly Alster: $\sone(\mathcal{G}_{K},\mathcal{G}_{{D}_{\Gamma}})$}
\end{center}

\begin{quote}
$\mathcal{G}_{D_{\Gamma}}$ is the family of infinite sets $\mathcal{U}$ where each member of $\mathcal{U}$ is a ${\sf G}_{\delta}$ subset of $X$, and for each nonempty open subset $U$ of $X$,  $\{V\in\mathcal{U}: U\cap V =\emptyset\}$ is finite.
\end{quote}

A space with the property $\sone(\mathcal{G}_K,\mathcal{G}_{\Gamma})$ has the property $\sone(\mathcal{G}_K,\mathcal{G}_{D_{\Gamma}})$. Indeed:
\begin{lemma}\label{denseetc}
If a space $X$ has a dense subset which has the property $\sone(\mathcal{G}_K,\mathcal{G}_{\Gamma})$, then $X$ has the property $\sone(\mathcal{G}_K,\mathcal{G}_{D_\Gamma})$.
\end{lemma}
\Proof Let $X$ be a space which has a dense subspace $Y$ which has the property $\sone(\mathcal{G}_K,\mathcal{G}_{\Gamma})$. If $(\mathcal{U}_n:n\in\naturals)$ is a sequence of elements of $\mathcal{G}_K$ for $X$, then choose for each $n\in \naturals$ an element $U_n\in\mathcal{U}$ such that  $\{U_n:n\in \naturals\}$ is an element of $\mathcal{G}_{{\Gamma}}$ for $Y$. Since $Y$ is dense in $X$, $\{U_n:n\in \naturals\}$ is an element of $\mathcal{G}_{D_{\Gamma}}$ for $X$. $\Box$

\begin{theorem}\label{weakprod} If $X$ has the property $\sone(\mathcal{G}_{K},\mathcal{G}_{{D}_{\Gamma}})$, then it is
\begin{enumerate}
  \item{productively weakly-Menger.}
  \item{productively weakly Hurewicz.}
\end{enumerate}
\end{theorem}
\Proof
We give the argument for weakly-Menger productive. Thus, let $Y$ be a weakly Menger space, and let $(\mathcal{U}_n:n\in\naturals)$ be a sequence of open covers of $X\times Y$. We may assume that each $\mathcal{U}_n$ is closed under finite unions.

For each $n$: For each compact subset $C$ of $X$ find for each $y\in Y$ a set $U_n(C,y)\in\mathcal{U}_n$ such that $C\times\{y\}\subseteq U_n(C,y)$. Then for each $n$ we find open subsets $V_n(C,y)\subset X$ and $W_n(C,y)\subset Y$ such that
\[
  C\times\{y\}\subseteq V_n(C,y)\times W_n(C,y)\subset U_n(C,y).
\]
For each $C$ define $\mathcal{U}_n(C) = \{W_n(C,y):y\in Y\}$, an open cover of $Y$.

Partition $\naturals$ into countably many infinite subsets $S_m$, $m\in\naturals\}$.

Since $Y$ is weakly Menger we find for each $m$, for each $n\in S_m$ a finite set $F(C,n)\subset Y$ such that $\bigcup\{W_n(C,y):y\in F(C,n),\, n\in S_m\}$ is dense in $Y$ in the following sense: For each nonempty open $T\subset Y$ there are infinitely many $n\in S_m$ for which $W_n(C,y)\cap T\neq \emptyset$ for some $y\in F(C,n)$.

Then for each $C$ define $V(C) =\bigcap_{n\in\naturals}\bigcap\{V_n(C,y):y\in F(C,n)\}$, a ${\sf G}_{\delta}$ subset of $X$ that contains $C$. Then $\mathcal{V} = \{V(C):C\subset X \mbox{ compact}\}$ is a member of $\mathcal{G}_K$ for $X$. Applying the fact that $X$ has the property $\sone(\mathcal{G}_K,\mathcal{G}_{D_{\Gamma}})$, we find a countable set $\{V(C_m):m\in\naturals\}$ such that for each nonempty open set $S\subset X$, for all but finitely many $n$, $S\cap V(C_m)\neq \emptyset$.

For each $n$ find the $m$ with $n\in S_m$ and define $\mathcal{V}_n = \{V_n(C_m,y)\times W_n(C_m,y) :y\in F(C_m,n)\}$. Then each $\mathcal{V}_n$ is a refinement of a finite subset of $\mathcal{U}_n$. We must show that $\bigcup_{n\in\naturals}(\bigcup\mathcal{V}_n)$ is dense in $X\times Y$. Thus, let $U\times V$ be a nonempty open subset of $X\times Y$. Choose $N$ so large that for all $m\ge N$ we have $V(C_m)\cap U\neq \emptyset$. For such an $m$ there are infinitely many $n\in S_m$ such that $W_n(C_m,y)\cap V\neq\emptyset$ for some $y\in F(C_m,n)$. Since $V(C_m)$ is a subset of $V_n(C_m)$, it follows that for infinitely many $n$ there are $y\in F(C_m,n)$ for which $V_n(C_m)\times W_n(C_m,y)\cap U\times V\neq \emptyset$.
$\Box$

As in the case of Lindel\"of productiveness, a proof of the following conjectures may depend on a characterization of the
\begin{conjecture}
If a space is productively Menger then it is productively weakly Menger.\\
If a space is productively Hurewicz then it is productively weakly Hurewicz.\\
\end{conjecture}

\begin{theorem}\label{weakrothbprod} If $X$ is a Rothberger space with the property $\sone(\mathcal{G}_K,\mathcal{G}_{D_{\Gamma}})$, then it is productively weakly Rothberger.
\end{theorem}
\Proof
Let $Y$ be a weakly Rothberger space. Let $(\mathcal{U}_n:n\in\naturals)$ be a sequence of open covers of $X\times Y$.
We may assume that for each $n$ the elements of $\mathcal{U}_n$ are of the form $U\times V$. Partition $\naturals$ into infinitely many infinite subsets  $S_m$, $m\in\naturals$. Then partition each $S_m$ into infinitely many infinite pairwise disjoint subsets $S_{m,k}$, $k\in\naturals$.

For each $m$, and for each compact subset $C$ of $X$, do the following. For each $k$ and for each $y\in Y$ choose $i_i<i_2<\cdots<i_t$ from $S_{m,k}$ such that for each $j$, $y\in V_{i_j}$, and such that $C\subseteq U_{i_1}\cup\cdots\cup U_{i_t}$. This is possible since $C$, a compact subset of $X$ is a Rothberger space. Define $V_{m,k}(y,C)$ to be the the intersection of these $V_{i_j}$, and $U_{m,k}(y,C)$ to be the union of these $U_{i_j}$. Then $\mathcal{W}_{m,k}(C)=\{V_{m,k}(y,C):y\in Y\}$ is an open cover of $Y$ for each $m$ and $k$. Applying the fact that $Y$ is weakly Rothberger to $(\mathcal{W}_{m,k}:k\in\naturals)$, choose for each $k$ a $y_k\in Y$ such that for each open set $V\subset Y$ there are infinitely many $k$ with $V_{m,k}(y_k,C)\cap V\neq\emptyset$. Define $V_m(C)$ to be the set $\bigcap_{k\in\naturals}U_{m,k}$.

Then for each compact set $C\subset X$ define $V(C) = \bigcap_{m\in\naturals}V_m(C)$, a ${\sf G}_{\delta}$ subset of $X$ which contains $C$. Now apply $\sone(\mathcal{G}_K,\mathcal{G}_{D_{\Gamma}})$ to the member $\{V(C):C\subset X\mbox{ xompact}\}$ of $\mathcal{G}_K$. We find a sequence $(V(C_m):m\in\naturals)$ of ${\sf G}_{\delta}$ sets such that for each nonempty open subset $U$ of $X$, for all but finitely many $n$, $V(C_n)\cap U$ is nonempty. Then the sequence $(V_n(C_n):n\in\naturals)$ has the same properties.

Now consider the sets $U_{m,k}(y_k,C_m)\times V_{m,k}(y_k,C_m)$ for $k,m\in\naturals$. For each nonempty open $U\times V\subset X\times Y$ there is an $N$ such that for all $m>N$, $V(C_m)\cap U$ is nonempty, whence $U_{m,k}\cap C$ is nonempty for all $k$ But for such an $m$, for infinitely many $k$ also $V_{m,k}(y_k,C_m)\cap V$ is nonempty. Thus, for infinitely many $m$ and $k$, $U_{m,k}(y_k,C_m)\times V_{m,k}(y_k,C_m)$ has nonempty intersection with $U\times V$. But now $U_{m,k}(y_k,C_m)$ is of the form $U_{i_1}\cup\cdots\cup U_{i_t}$ while $V_{m,k}(y_k,C_m)$ is of the form $V_{i_1}\cap\cdots\cap V_{i_t}$, where $i_1<\cdots<i_t$ are from $S_{m,k}$ and $U_{i_j}\times V_{i_j}$ is an element of $\mathcal{U}_{i_j}$.

It follows that there is a sequence of $S_n\in\mathcal{U}_n$ such that for each nonempty open $U\times V\subset X\times Y$, here are infinitely many $n$ with $S_n\cap U\times V$ nonempty. $\Box$

\begin{corollary}\label{sigmcptwkalster} Every space which has a dense $\sigma$-compact subset has property $\sone(\mathcal{G}_K,\mathcal{G}_{{D}_{\Gamma}})$.
\end{corollary}
\Proof A $\sigma$-compact space has the property $\sone(\mathcal{G}_K,\mathcal{G}_{{\Gamma}})$. $\Box$

\begin{corollary}\label{sepwkalster} Every separable space has the property $\sone(\mathcal{G}_K,\mathcal{G}_{{D}_{\Gamma}})$.
\end{corollary}
\Proof Let $X$ be a separable space, and let $D$ be a countable dense subset of $X$. Then $D$ is a dense $\sigma$-compact subset of $X$. $\Box$

\begin{corollary}\label{sepprod}
Every separable space is productively weakly Lindel\"of, productively weakly Menger and productively weakly Hurewicz.
\end{corollary}
\Proof Theorem \ref{weakprod} and Corollary \ref{sepwkalster}. $\Box$

\begin{corollary}\label{Rpowers} For each cardinal number $\kappa$, $\reals^{\kappa}$ is productively weakly Lindel\"of, productively weakly Menger, productively weakly Hurewicz and productively weakly Rothberger.
\end{corollary}
\Proof We only consider the case when $\kappa$ in infinite. By Proposition 4 of \cite{Corson},
\[
  \{f\in\reals^{\kappa}: \vert\{i:f(i)\neq 0\}\vert<\aleph_0\}
\]
is a dense $\sigma$-compact subset of $\reals^{\kappa}$. By Lemma  \ref{weakalsterLindelof}, Corollary \ref{sigmcptwkalster} and  Theorem \ref{weakprod} $\reals^{\kappa}$ has the claimed properties.

Since $\rationals$ is $\sigma$-compact and Rothberger, Proposition 4 of \cite{Corson} implies that the subset
\[
  S_{\kappa} = \{f\in ^{\kappa}\rationals:\, \vert\{\alpha\in\kappa:\, f(\alpha)=0\}\vert <\aleph_0\}
\]
of $\rationals^{\kappa}$ is $\sigma$-compact. By Corollary 15 of \cite{MSRBGs}, TWO has a winning strategy in the game $\gone(\Omega,\Gamma)$ in $S_{\kappa}$, and thus $S_{\kappa}$ is a $\gamma$-space. It is evident that $S_{\kappa}$ is a dense subset of $\reals^{\kappa}$.

But $S_{\kappa}$ is dense in $\reals^{\kappa}$ and so as $S_{\kappa}$ is productively weakly Rothberger by Theorem \ref{weakrothbprod}, also $\reals^{\kappa}$ is productively weakly Rothberger. $\Box$

Separable spaces are weakly Alster, and in fact have several additional properties: For example a separable space is also weakly Rothberger. Note that a compact space could be weakly Rothberger without being Rothberger: The closed unit interval is a compact separable space but is not a Rothberger space.

\begin{corollary}\label{weaklyproductive}
If $X$ is a separable space then $X$ is productively weakly Rothberger.
\end{corollary}
\Proof Every countable space has the property $\sone(\mathcal{G}_K,\mathcal{G}_{\Gamma})$, and thus the property $\sone(\mathcal{G}_K,\mathcal{G}_{D_{\Gamma}})$. Since a countable space is Rothberger, Theorem \ref{weakrothbprod} implies that any countable space is productively weakly Rothberger. Apply Theorem \ref{denseupwabsolute} to conclude that any separable space is productively weakly Rothberger. $\Box$

\begin{lemma}\label{cptrothb} Compact Rothberger spaces are productively weakly Rothberger.
\end{lemma}
\Proof Let $X$ be a compact Rothberger space and let $Y$ by weakly Rothberger. Let $(\mathcal{U}_n:n\in\naturals)$ be a sequence of open covers of $X\times Y$. We may assume each $\mathcal{U}_n$ consist of sets of the form $U\times V$. Write $\naturals = \bigcup_{n\in\naturals}S_n$ where the $S_n$'s are infinite and pairwise disjoint.

Fix $n$ as well as $y\in Y$, and choose a finite sequence $U_{i_1},\cdots,U_{i_k}$ where
\begin{enumerate}
  \item{$i_1<\cdots< i_k$ are elements of $S_n$;}
  \item{$U_{i_j}$ is an element of $\mathcal{U}_{i_j}$, $1\le j\le k$;}
  \item{$X =  U_{i_1}\cup\cdots\cup U_{i_k}$;}
  \item{$y\in V_{i_1}\cap\cdots\cap V_{i_k} = V_n(y)$, say.}
\end{enumerate}
This defines for each $n$ an open cover $\mathcal{H}_n = \{V_n(y):y\in Y\}$ of $Y$. Now apply the fact that $Y$ is weakly Rothberger to the sequence $(\mathcal{H}_n:n\in\naturals)$ and select for each $n$ an $H_n\in \mathcal{H}_n$ such that for each open set $V\subseteq Y$, for infinitely many $n$, $V\cap H_n$ is nonempty. This produces a sequence of elements $U_n\in\mathcal{U}_n$, $n\in\naturals$, such that $\bigcup\{U_n:n\in\naturals\}$ is dense in $X\times Y$. $\Box$

\begin{problem}\label{wkRothbProductive} Is it true that if $X$ has the property $\sone(\mathcal{G}_K,\mathcal{G}_{D_{\Gamma}})$ and each compact subset of $X$ is weakly Rothberger, then $X$ is productively weakly Rothberger?
\end{problem}

Using techniques from the previous section we can prove:
\begin{lemma}\label{weakrothb} If each compact subspace of $X$ is a weakly Rothberger space and if $X$ satisfies $\sone(\open_K,\dense)$, then $X$ it is weakly Rothberger.
\end{lemma}

A space is \emph{weakly Gerlits-Nagy space} if it satisfies $\sone(\Omega,\dense^{gp})$.

\begin{lemma} $\sone(\Omega, \dense^{gp})$ is equivalent to $\sone(\Lambda, \dense^{gp})$.
\end{lemma}
\Proof Since $\Omega\subset \Lambda$ we need to prove that is $X$ has $\sone(\Omega, \dense^{gp})$ then is has $\sone(\Lambda, \dense^{gp})$. Let $({\mathcal U}_n:n\in\naturals)$ be  a sequence of $\lambda$-covers. Then
${\mathcal V}_n=\{U\in {\mathcal U}_n: F\subset U$, $F$ finite subset of $X\}$.
Then $({\mathcal V}_n:n\in\naturals)$ is a sequence of $\omega$-cover of $X$. Apply the hypothesis we can select an element $U_n\in {\mathcal V}_n$ such that $\{U_n:n\in\naturals\}$ is an element of $\dense^{gp}$.
Then since $U_n$ is extract for each ${\mathcal U}_n$ then $\{U_n:n\in\naturals\}$ witnesses that $X$ satisfies $\sone(\Lambda, \dense^{gp})$. $\Box$

\begin{lemma} \label{cptGN} Compact Gerlits-Nagy spaces are productively weakly Gerlits-Nagy.
\end{lemma}

\Proof Let $X$ be a compact Gerlits-Nagy space and let $Y$ be a weakly Gerlits-Nagy. Let $(\mathcal{U}_n:n\in\naturals)$ be a sequence of $\omega$-covers of $X\times Y$. We may assume each $\mathcal{U}_n$ consist of sets of the form $U\times V$. Write $\naturals = \bigcup_{n\in\naturals}S_n$ where the $S_n$'s are infinite and pairwise disjoint.

Fix $n$ as well as $F$ finite subset of $Y$, and choose a finite sequence $U_{i_1},\cdots,U_{i_k}$ where
\begin{enumerate}
  \item{$i_1<\cdots< i_k$ are elements of $S_n$;}
  \item{$U_{i_j}$ is an element of $\mathcal{U}_{i_j}$, $1\le j\le k$;}
  \item{$X =  U_{i_1}\cup\cdots\cup U_{i_k}$;}
  \item{$F\subseteq  V_{i_1}\cap\cdots\cap V_{i_k} = V_n(F)$, say.}
\end{enumerate}
This defines for each $n$ an $\omega$-cover $\mathcal{H}_n = \{V_n(F):F$ finite subset of $Y\}$ of $Y$. Now apply the fact that $Y$ is weakly Gerlits-Nagy to the sequence $(\mathcal{H}_n:n\in\naturals)$ and select for each $n$ an $H_n\in \mathcal{H}_n$ such that for each open set $V\subseteq Y$, for infinitely many $n$, $V\cap H_n$ is nonempty. This produces a sequence of elements $U_n\in\mathcal{U}_n$, $n\in\naturals$, such that $\bigcup\{U_n:n\in\naturals\}$ is dense in $X\times Y$. $\Box$

\begin{problem} Is each compact weakly Rothberger space a weakly Gerlits-Nagy space? (or maybe a weakly $\gamma$-space?)
\end{problem}

We also don't know if weakly Lindel\"of {\sf P}-spaces have stronger properties in terms of products. They are at least productively weakly Lindel\"of as we now show:

\begin{lemma}\label{weakLandP}
If $X$ is a weakly Lindel\"of {\sf P}-space, then it has property $\sone(\mathcal{G}_K,\mathcal{G}_{D})$.
\end{lemma}
\Proof For each $n\in\naturals$ let a set $\mathcal{U}_n\in\mathcal{G}_K$ be given. Since $X$ is a {\sf P}-space, $\mathcal{U}_n$ is an open cover of $X$. For each compact subset $C$ of $X$ choose for each $n$ a $U_n(C)\in\mathcal{U}_n$ with $C\subseteq U_n(C)$, and then define $V(C) = \bigcap_{n\in\naturals}U_n(C)$. Then $V(C)$ is a ${\sf G}_{\delta}$ subset of $X$ and $\mathcal{V} = \{V(C):C\subset X \mbox{ compact}\}$ is an open cover of $X$ as $X$ is a {\sf P}-space. Since $X$ is weakly Lindel\"of, choose a countable subset $\{V(C_n):n\in\naturals\}$ of $\mathcal{U}$ which has the property that there is for each nonempty open $V\subset X$ an infinite number of $n$ with $V\cap V(C_n)\neq\emptyset$. 
For each $n$, choose $U_n = U_n(C_n)\in\mathcal{U}_n$. Then $\{U_n:n\in\naturals\}$ is a member of $\mathcal{G}_D$.
$\Box$

It also follows that weakly Lindel\"of {\sf P}-spaces are weakly Rothberger. Since the finite product of {\sf P}-spaces are {\sf P}-spaces, Lemma \ref{weakLandP} and Lemma \ref{weakalsterLindelof} imply that the product of finitely many weakly Lindel\"of {\sf P}-spaces is a weakly Lindel\"of {\sf P}-space.

\begin{problem} Let $X$ be a weakly Lindel\"of {\sf P}-space.
\begin{itemize}
  \item[(a)]{Does $X$ then have the property $\sone(\mathcal{G}_K,\mathcal{G}_{D_{\Gamma}})$?}
  \item[(b)]{Is $X$ productively weakly-Rothberger?}
  \item[(c)]{Is $X$ weakly-Hurewicz?}
\end{itemize}
\end{problem}

For the family $\mathcal{D}$ we introduce
\[
  \dense_{\Gamma}=\{\mathcal{U}\in\dense:\, \mathcal{U} \mbox{ infinite and each infinite subset of $\mathcal{U}$ is in }\dense\}.
\]

A space is \emph{weak $\gamma$-space} if it satisfies $\sone(\Omega,\dense_{\Gamma})$.

\begin{lemma}\label{cptgamma} Compact $\gamma$-spaces are productively weakly $\gamma$-spaces.
\end{lemma}
\Proof Let $X$ be a compact $\gamma$-pace and let $Y$ be a weakly $\gamma$-space. Let $(\mathcal{U}_n:n\in\naturals)$ be a sequence of $\omega$-covers of $X\times Y$. We may assume each $\mathcal{U}_n$ consist of sets of the form $U\times V$. Write $\naturals = \bigcup_{n\in\naturals}S_n$ where the $S_n$'s are infinite and pairwise disjoint.

Fix $n$ as well as $F$ finite subset of $Y$, and choose a finite sequence $U_{i_1},\cdots,U_{i_k}$ where
\begin{enumerate}
  \item{$i_1<\cdots< i_k$ are elements of $S_n$;}
  \item{$U_{i_j}$ is an element of $\mathcal{U}_{i_j}$, $1\le j\le k$;}
  \item{$X =  U_{i_1}\cup\cdots\cup U_{i_k}$;}
  \item{$F\subseteq  V_{i_1}\cap\cdots\cap V_{i_k} = V_n(F)$, say.}
\end{enumerate}
This defines for each $n$ an $\omega$-cover $\mathcal{H}_n = \{V_n(F):F$ finite subset of $Y\}$ of $Y$. Now apply the fact that $Y$ is weakly $\gamma$-space to the sequence $(\mathcal{H}_n:n\in\naturals)$ and select for each $n$ an $H_n\in \mathcal{H}_n$ such that $\{H_n:n\in\naturals\}\in {\mathcal D}_\Gamma$. This produces a sequence of elements $U_n\in\mathcal{U}_n$, $n\in\naturals$, such that $\bigcup\{U_n:n\in\naturals\}$ is dense in $X\times Y$. $\Box$

\section*{Acknowledgements}

Part of this work was done while B. Pansera was visiting the Department of Mathematics at Boise State University. The Department's hospitality and support during this visit is gratefully acknowledged.


\end{document}